\documentstyle[12pt,amssymb,amstex]{amsart}

\numberwithin{equation}{section}

\def\ga{{\frak A}}

\def\bc{{\mathbb C}}

\def\bn{{\mathbb N}}

\def\a{\alpha}

\def\e{\epsilon}

\def\p{\psi}

\def\f{\varphi}

 \def\Om{\Omega}

\def\xb{{\mathbf{x}}}
\def\yb{{\mathbf{y}}}

\newtheorem{thm}{Theorem}[section]
\newtheorem{lem}[thm]{Lemma}

\newtheorem{defin}[thm]{Definition}

\def\id{\mathop{\rm id}}

\def\id{{\bf 1}\!\!{\rm I}}

\begin{document}

\small

\title[On noncommutative unique ergodicity]
{A note on noncommutative unique ergodicity\\ and weighted means}

\author{Luigi Accardi}
\address{Luigi Accardi\\
Centro Interdisciplinare Vito Volterra\\
II Universit\`{a} di Roma ``Tor Vergata''\\
Via Columbia 2, 00133 Roma, Italy} \email{{\tt
accardi@@volterra.uniroma2.it}}
\author{Farrukh Mukhamedov}
\address{Farrukh Mukhamedov\\
 Department of Computational \& Theoretical Sciences\\
Faculty of Sciences, International Islamic University Malaysia\\
P.O. Box, 141, 25710, Kuantan\\
Pahang, Malaysia}\email{{\tt far75m@@yandex.ru},{\tt
farrukh\_m@@iiu.edu.my}}

\begin{abstract}
In this paper we study unique ergodicity of $C^*$-dynamical system
$(\ga,T)$, consisting of a unital $C^*$-algebra $\ga$ and a Markov
operator $T:\ga\mapsto\ga$, relative to its fixed point subspace, in
terms of Riesz summation which is weaker than Cesaro one. Namely, it
is proven that  $(\ga,T)$ is uniquely ergodic relative to its fixed
point subspace if and only if its Riesz means
\begin{equation*}
\frac{1}{p_1+\cdots+p_n}\sum_{k=1}^{n}p_kT^kx
\end{equation*}
converge to $E_T(x)$ in $\ga$ for any $x\in\ga$, as $n\to\infty$,
here $E_T$ is an projection of $\ga$ to the fixed point subspace of
$T$. It is also constructed a uniquely ergodic entangled Markov
operator relative to its fixed point subspace, which is not ergodic.
 \vskip 0.3cm \noindent {\it
Mathematics Subject
Classification}: 47A35, 46L35, 46L55.\\
{\it Key words}: uniquely ergodic, Markov operator, Riesz means.
\end{abstract}

\maketitle

\footnotetext[1]{The author (F.M.) is on leave from National
University of Uzbekistan, Department of Mechanics \& Mathematics,
Vuzgorodok, 100174 Tashkent,  Uzbekistan }

\section{Introduction}

It is known \cite{KSF, W} that one of the important notions in
ergodic theory is unique ergodicity of a homeomorphism $T$ of a
compact Hausdorff space $\Om$. Recall that $T$ is {\it uniquely
ergodic} if there is a unique $T$--invariant Borel probability
measure $\mu$ on $\Omega$. The well known Krylov-Bogolyubov theorem
\cite{KSF} states that $T$ is uniquely ergodic if and only if for
every $f\in C(\Om)$ the averages
\begin{equation*}
\frac{1}{n}\sum_{k=0}^{n-1}f(T^kx)
\end{equation*}
converge uniformly to the constant $\int f\,d\mu$, as
$n\to\infty$.

The study of ergodic theorems in recent years showed that the
ordinary Cesaro means have been replaced by weighted averages
\begin{equation}\label{a_k}
\sum_{k=0}^{n-1}a_kf(T^kx).
\end{equation}
Therefore, it is natural to ask: is there a weaker summation than
Cesaro, ensuring the unique ergodicity. In \cite{K} it has been
established that unique ergodicity implies uniform convergence of
\eqref{a_k}, when $\{a_k\}$ is Riesz weight (see also \cite{I} for
similar results). In \cite{BLRT} similar problems were considered
for transformations of Hilbert spaces.

On the other hand, since the theory of quantum dynamical systems
provides a convenient mathematical description of irreversible
dynamics of an open quantum system (see \cite{AH},\cite{BR})
investigation of ergodic properties of such dynamical systems have
had a considerable growth. In a quantum setting, the matter is more
complicated than in the classical case. Some differences between
classical and quantum situations are pointed out in
\cite{AH},\cite{NSZ}. This motivates an interest to study dynamics
of quantum systems (see \cite{FR1,FR2,FV}). Therefore, it is then
natural to address the study of the possible generalizations to
quantum case of various ergodic properties known for classical
dynamical systems. In \cite{LP},\cite{MT} a non-commutative notion
of unique ergodicity was defined, and certain properties were
studied. Recently in \cite{AD} a general notion of unique ergodicity
for automorphisms of a $C^*$-algebra with respect to its fixed point
subalgebra has been introduced. The present paper is devoted to a
generalization of such a notion for positive mappings of
$C^*$-algebras, and its characterization in term of Riesz means.

The paper is organized as follows: section 2 is devoted to
preliminaries, where we recall some facts about $C^*$-dynamical
systems and the Riesz summation of a sequence on $C^*$-algebras.
Here we define a notion of unique ergodicity of $C^*$-dynamical
system relative to its fixed point subspace. In section 3 we prove
that a $C^*$-dynamical system $(\ga,T)$ is uniquely ergodic relative
to its fixed point subspace if and only if its Riesz means (see
below)
\begin{equation*}
\frac{1}{p_1+\cdots+p_n}\sum_{k=1}^{n}p_kT^kx
\end{equation*}
converge to $E_T(x)$ in $\ga$ for any $x\in\ga$, here $E_T$ is a
projection of $\ga$ onto the fixed point subspace of $T$. Note
however that if $T$ is completely positive then $E_T$ is a
conditional expectation (see \cite{C1,P}. On the other hand it is
known \cite{MT} that unique ergodicity implies ergodicity.
Therefore, one can ask: can a $C^*$-dynamical system which is
uniquely ergodic relative to its fixed point subspace  be ergodic?
It turns out that this question has a negative answer. More
precisely, in section 4 we construct entangled Markov operator which
is uniquely ergodic  relative to its fixed point subspace, but which
is not ergodic.

\section{Preliminaries}

In this section we recall some preliminaries concerning
$C^*$-dynamical systems.

Let $\ga$ be a $C^*$-algebra with unit $\id$. An element $x\in\ga$
is called {\it positive} if there is an element $y\in\ga$ such
that $x=y^*y$. The set of all  positive elements will be denoted
by $\ga_+$. By $\ga^*$ we denote the conjugate space to $\ga$. A
linear functional $\f\in\ga^*$ is called {\it Hermitian} if
$\f(x^*)=\overline{\f(x)}$ for every $x\in\ga$. A Hermitian
functional $\f$ is called {\it state} if $\f(x^*x)\geq 0$ for
every $x\in\ga$ and $\f(\id)=1$. By $S_\ga$ (resp. $\ga^*_h$) we
denote the set of all states (resp. Hermitian functionals) on
$\ga$.  By $M_n(\ga)$ we denote the set of all $n\times
n$-matrices $a=(a_{ij})$ with entries $a_{ij}$ in $\ga$.
\begin{defin} A linear operator $T:\ga\mapsto\ga$ is called:
\begin{enumerate}
    \item[(i)] {\it positive}, if $Tx\geq 0$ whenever $x\geq 0$;
    \item[(ii)] {\it $n$-positive} if the linear mapping $T_n:M_n(\ga)\mapsto M_n(\ga)$ given by
$T_n(a_{ij})=(T(a_{ij}))$ is positive;
    \item[(iii)] {\it completely positive} if it is $n$-positive for all $n\in\bn$.
    \end{enumerate}
\end{defin}
A positive mapping $T$ with $T\id=\id$ is called {\it Markov
operator}. A pair $(\ga,T)$ consisting of a $C^*$-algebra $\ga$ and
a Markov operator $T:\ga\mapsto\ga$ is called {\it a $C^*$-dynamical
system}. The $C^*$-dynamical system $(\ga,\f,T)$ is called {\it
uniquely ergodic} if there is a unique invariant state $\f$ (i.e.
$\f(Tx)=\f(x)$ for all $x\in \ga$) with respect to $T$. Denote
\begin{equation}\label{fix}
\ga^T=\{x\in\ga: Tx=x\}.
\end{equation}
It is clear that $\ga^T$ is a closed linear subspace of $\ga$, but
in general it is not a subalgebra of $\ga$ (see sec. 3). We say that
$(\ga,T)$ is {\em uniquely ergodic relative to $\ga^T$} if every
state of $\ga^T$ has a unique $T$--invariant state extension to
$\ga$. In the case when $\ga^T$ consists only of scalar multiples of
the identity element, this reduces to the usual notion of unique
ergodicity. Note that for an automorphism such a notion has been
introduced in \cite{AD}.

Now suppose we are given a sequence of numbers $\{p_n\}$ such that
$p_1>0$, $p_k\geq 0$ with $\sum\limits_{k=1}^\infty p_k=\infty$. We
say that a sequence $\{s_n\}\subset\ga$ is  {\it Riesz convergent}
to an element $s\in\ga$ if the sequence
$$
\frac{1}{\sum_{k=1}^np_k}\sum_{k=1}^np_ks_k
$$
converges to $s$ in $\ga$, and it is denoted by $s_n\to s \
(R,p_n)$. The numbers $p_n$ are called {\it weights}. If $s_n\to
s$ implies $s_n\to s (R,p_n)$ then Riesz-converges is said to be
{\it regular}. The regularity condition (see \cite{H}, Theorem 14)
is equivalent to
\begin{equation}\label{reg}
\frac{p_n}{p_1+p_2+\cdots+p_n}\to 0  \qquad \textrm{as} \ \
n\to\infty.
\end{equation}

Basics about $(R,p_n)$ convergence can be found in \cite{H}.

Recall the following lemma which shows that Riesz convergence is
weaker than Cesaro convergence (see \cite{H},\cite{K}).

\begin{lem}\label{cr} (\cite{H}, Theorem 16) Assume that $p_{n+1}\leq p_n$ and
\begin{equation}\label{pp0}
\frac{np_n}{p_1+\cdots+p_n}\leq C \qquad \forall n\in\bn
\end{equation}
for some constant $C>0$. Then Cesaro convergence implies $(R,p_n)$
convergence. \end{lem}

\section{Unique ergodicity}

In this section we are going to characterize unique ergodicity
relative to $\ga^T$ of $C^*$-dynamical systems.  To do it we need
the following

\begin{lem}\label{inv} (cf. \cite{MT},\cite{AD})
Let $(\ga,T)$ be uniquely ergodic relative to $\ga^T$. If
$h\in\ga^*$ is invariant with respect to $T$ and
$h\upharpoonright\ga^T=0$, then  $h=0.$
\end{lem}
\begin{pf} Let us first assume that $h$ is Hermitian. Then there is a unique
Jordan decomposition \cite{T} of $h$ such that
\begin{equation}\label{her}
h=h_+-h_-, \ \ \|h\|_1=\|h_+\|_1+\|h_-\|_1,
\end{equation}
where $\|\cdot\|_1$ is the dual norm on $\ga^*$. The invariance of
$h$ implies that
$$
h\circ T=h_+\circ T-h_-\circ T=h_+-h_-.
$$
Using $\|h_+\circ T\|_1=h_+(\id)=\|h_+\|_1$, similarly $\|h_+\circ
T\|_1=\|h_+\|_1$, from uniqueness of the decomposition we find
$h_+\circ T=h_+$ and $h_-\circ T=h_-$. From $h\upharpoonright
\ga^T=0$ one gets $h(\id)=0$, which implies that
$\|h_+\|_1=\|h_-\|_1$. On the other hand, we also have
$\frac{h_+}{\|h_+\|_1}=\frac{h_-}{\|h_-\|_1}$ on $\ga^T$. So,
according to the unique ergodicity relative to $\ga^T$ we obtain
$h_+=h_-$ on $\ga$. Consequently,  $h=0$. Now let $h$ be an
arbitrary bounded, linear functional. Then it can be written as
$h=h_1+ih_2$, where $h_1$ and $h_2$ are Hermitian. Again invariance
of $h$ implies that $h_i\circ T=h_i$, $i=1,2$. From
$h\upharpoonright\ga^T=0$ one gets $h_k\upharpoonright\ga^T=0$,
$k=1,2$.  Consequently, according to the above argument, we obtain
$h=0$.
\end{pf}

Now we are ready to formulate a criterion for unique ergodicity of
$C^*$-dynamical system in terms of $(R,p_n)$ convergence. In the
proof we will follow some ideas used in \cite{AD, K, MT}.

\begin{thm}\label{unique}  Let $(\ga,T)$ be a $C^*$-dynamical system. Assume that the weight $\{p_n\}$ satisfies
\begin{equation}\label{pp}
P(n):=\frac{p_1+|p_2-p_1|+\cdots+|p_n-p_{n-1}|+p_n}{p_1+p_2+\cdots+p_n}\to
0 \ \ \textrm{as} \ \  n\to\infty.
\end{equation}
Then the following conditions are equivalent:
\begin{itemize}
\item[(i)] $(\ga,T)$ is uniquely ergodic relative to $\ga^T$;
\item[(ii)] The set $\ga^T+\{a-T(a):a\in\ga\}$ is
dense in $\ga$;
 \item[(iii)] For all $x\in \ga$,
$$
T^nx\to E_T(x) \ (R,p_n),
$$
where $E_T(x)$ is a positive norm one projection onto $\ga^T$ such
that $E_TT=TE_T=E_T$; Moreover, the following estimation holds:
\begin{equation}\label{estim}
\bigg\|\frac{1}{\sum_{k=1}^np_k}\sum_{k=1}^{n}p_kT^k(x)-E_T(x)\bigg\|\leq
P(n)\|x\|, \ \ \ n\in\bn,
\end{equation}
for every $x\in\ga$;
 \item[(iv)] For every $x\in\ga$ and $\p\in
S_\ga$
$$
\p(T^k(x))\to \p(E_T(x)) \ (R,p_n).
$$
\end{itemize}
\end{thm}

\begin{pf} Consider the implication (i)$\implies$(ii). Assume that
$\overline{\ga^T+\{a-T(a):a\in\ga\}}\neq\ga$; then there is an
element $x_0\in\ga$ such that
$x_0\notin\overline{\ga^T+\{a-T(a):a\in\ga\}}$. Then according to
the Hahn-Banach theorem there is a functional $h\in\ga^*$ such that
$h(x_0)=1$ and
$h\upharpoonright\overline{\ga^T+\{a-T(a):a\in\ga\}}=0$. The last
condition implies that $h\upharpoonright\ga^T=0$ and $h\circ T=h$.
Hence, Lemma \ref{inv} yields that $h=0$, which contradicts to
$h(x_0)=1$.

(ii)$\implies$(iii): It is clear that for every element of the form
$y=x-T(x)$, $x\in\ga$ by \eqref{pp} we have
\begin{eqnarray}\label{lim1}
\frac{1}{\sum_{k=1}^np_k}\bigg\|\sum_{k=1}^{n}p_kT^k(y)\bigg\| &=&
\frac{1}{\sum_{k=1}^np_k}\bigg\|\sum_{k=1}^np_k(T^{k+1}(x)-T^kx)\bigg\|\nonumber\\
&=&
\frac{1}{\sum_{k=1}^np_k}\|p_1Tx+(p_2-p_1)T^2x+\cdots\nonumber\\
&&+(p_n-p_{n-1})T^nx-p_nT^{n+1}x\|\nonumber\\
&\leq& P(n)\|x\| \to 0 \ \ \textrm{as} \ \ n\to\infty.
\end{eqnarray}

 Now let $x\in\ga^T$, then
\begin{equation}\label{lim2}
\lim_{n\to\infty}\frac{1}{\sum_{k=1}^np_k}\sum_{k=1}^{n}p_kT^k(x)=x.
\end{equation}

Hence, for every $x\in \ga^T+\{a-T(a):a\in\ga\}$ the limit
$$
\lim_{n\to\infty}\frac{1}{\sum_{k=1}^np_k}\sum_{k=1}^{n}p_kT^kx
$$
exists, which is denoted by $E_T(x)$. It is clear that $E_T$ is a
positive linear operator from  $\ga^T+\{a-T(a):a\in\ga\}$ onto
$\ga^T$. Positivity and $E_T\id=\id$ imply that $E_T$ is bounded.
From \eqref{lim1} one  obviously gets that $E_TT=TE_T=E_T$.
According to (ii) the operator $E_T$ can be uniquely extended to
$\ga$, this extension is denoted by the same symbol $E_T$. It is
evident that $E_T$ is a positive projection with  $\|E_T\|=1$.

Now take an arbitrary $x\in\ga$. Then again using (ii), for any
$\e>0$ we can find $x_\e\in \ga^T+\{a-T(a):a\in\ga\}$ such that
$\|x-x_\e\|\leq\e.$ By means of \eqref{lim1},\eqref{lim2} we
conclude that
$$
\bigg\|\frac{1}{\sum_{k=1}^np_k}\sum\limits_{k=1}^{n}p_kT^k(x_\e)-E_T(x_\e)\bigg\|\leq
P(n)\|x_\e\|
$$
Hence, one has
\begin{eqnarray*}
\bigg\|\frac{1}{{\sum_{k=1}^np_k}}\sum_{k=1}^{n}p_kT^k(x)-E_T(x)\bigg\|&\leq
&
\bigg\|\frac{1}{\sum_{k=1}^np_k}\sum_{k=1}^{n}p_kT^k(x-x_\e)\bigg\|\\
&&+\bigg\|\frac{1}{\sum_{k=1}^np_k}\sum_{k=1}^{n}p_kT^k(x_\e)-E_T(x_\e)\bigg\|\\
&&+\|E_T(x-x_\e)\|\\
&\leq &2\|x-x_\e\|+P(n)\|x_\e\|\\
&\leq &P(n)\|x\|+(2+P(n))\e
\end{eqnarray*}
which with the arbitrariness of $\e$ implies \eqref{estim}.

 Consequently,
$$
\lim_{n\to\infty}\frac{1}{\sum_{k=1}^np_k}\sum_{k=1}^{n}p_kT^kx=E_T(x)
$$
is valid for every $x\in\ga$.

The mapping $E_T$ is a unique $T$-invariant positive projection.
Indeed, if $\tilde E:\ga\to\ga^T$ is any $T$-invariant positive
projection onto $\ga^T$, then
\begin{equation*}
\tilde E(x)=\frac{1}{\sum_{k=1}^np_k}\sum_{k=1}^{n}p_k\tilde
E(T^k(x))=\tilde
E\bigg(\frac{1}{\sum_{k=1}^np_k}\sum_{k=1}^{n}p_kT^k(x)\bigg).
\end{equation*}
Taking the limit as $n\to\infty$ gives
\begin{equation*}
\tilde E(x)=\tilde E(E_T(x))=E_T(x).
\end{equation*}

The implication (iii)$\implies$(iv) is obvious. Let us consider
(iv)$\implies$(i). Let $\p$  be any state on $\ga^T$,  then $\p\circ
E_T$ is a $T$-invariant extension of $\p$ to $\ga$. Assume that
$\phi$ is any $T$--invariant, linear extension of $\p$. Then
\begin{equation*}
\phi(x)=\frac{1}{\sum_{k=1}^np_k}\sum_{k=1}^{n}p_k\phi(T^k(x))=\phi\bigg(\frac{1}{\sum_{k=1}^np_k}\sum_{k=1}^{n}p_kT^k(x)\bigg).
\end{equation*}
Now taking the limit from both sides of the last equality as
$n\to\infty$ one gives
\begin{equation*}
\phi(x)=\phi(E_T(x))=\p(E_T(x)),
\end{equation*}
so $\phi=\p\circ E_T$.
\end{pf}

{\bf Remark 1.} If we choose $p_n=1$ for all $n\in\bn$ then it is
clear that the condition \eqref{pp} is satisfied, hence we infer
that unique ergodicity relative to $\ga^T$ is equivalent to the
norm convergence of the mean averages, i.e.
$$
\frac{1}{n}\sum_{k=1}^nT^k(x),
$$
which recovers the result of \cite{AD}.

{\bf Remark 2.} If the condition \eqref{pp0} is satisfied then
condition \eqref{pp} is valid as well. This means that unique
ergodicity would remain true if Cesaro summation is replaced by a
weaker. Theorem \ref{unique} extends a result of \cite{MT}.

{\bf Example. } If we define $p_n=n^\a$ with $\a>0$, then one can
see that $\{p_n\}$ is an increasing sequence and condition
\eqref{pp} is also satisfied. This provides a concrete example of
weights.

{\bf Remark 3.} Note that some nontrivial examples of uniquely
ergodic quantum dynamical systems based on automorphisms, has been
given in \cite{AD}. Namely, it was proved that free shifts based on
reduced $C^{*}$--algebras of RD--groups (including the free group on
infinitely many generators), and amalgamated free product
$C^{*}$--algebras, are uniquely ergodic relative to the fixed--point
subalgebra. In \cite{FM} it has been proved that such shifts possess
a stronger property called $F$-strict weak mixing (see also
\cite{MT}).

{\bf Observation.} We note that, in general, the projection $E_T$ is
not a conditional expectation, but when $T$ is an automorphism then
it is so.  Now we are going to provide an example of Markov operator
which is uniquely ergodic  relative to its fixed point subspace for
which the projector $E_T$ is not a conditional expectation.

Consider the algebra $M_d(\bc)$ - $d\times d$ matrices over $\bc$.
For a matrix $\xb=(x_{ij})$ by $\xb^t$ we denote its transpose
matrix, i.e. $\xb^t=(x_{ji})$. Define a mapping $\phi:M_d(\bc)\to
M_d(\bc)$ by $\phi(\xb)=\xb^t$. Then it is known \cite{P} that such
a mapping is positive, but not completely positive. One can see that
$\phi$ is a Markov operator. Due to the equality
$$
\xb=\frac{\xb+\xb^t}{2}+\frac{\xb-\xb^t}{2}
$$
condition (ii) of Theorem \ref{unique} is satisfied, so $\phi$ is
uniquely ergodic with respect to $M_d(\bc)^\phi$. Hence, the
corresponding projection $E_\phi$ is given by
$E_\phi(\xb)=(\xb+\xb^t)/2$, which is not completely positive.
Moreover,  $M_d(\bc)^\phi$ is the set of all symmetric matrices,
which do not form an algebra. So, $E_\phi$ is not a conditional
expectation.

\section{A uniquely ergodic entangled Markov operator}

In recent developments of quantum information many people have
discussed the problem of finding a satisfactory quantum
generalization of classical random walks. Motivating this in
\cite{AF, F} a new class of quantum Markov chains was constructed
which are at the same time purely generated and uniquely determined
by a corresponding classical Markov chain. Such a class of Markov
chains was constructed by means of entangled Markov operators. In
one's turn they were associated with Schur multiplication. In that
paper, ergodicity and weak clustering properties of such chains were
established. In this section we are going to provide entangled
Markov operator which is uniquely ergodic  relative to its fixed
point subspace, but which is not ergodic.

Let us recall some notations. To define Schur multiplication, we
choose an orthonormal basis $\{e_j\}$, $j=1,...,d$ in a
$d$-dimensional Hilbert space $H_d$ which is kept fixed during the
analysis. In such a way, we have the natural identification $H_d$
with  $C^d$. The corresponding system of matrix units $e_{ij} =
e_i\otimes e_j$ identifies $B(H_d)$ with $M_d(\bc)$. Then, for $\xb
= \sum_{i,j=1}^d x_{ij}e_{ij}$, $\yb =\sum_{i,j=1}^d y_{ij}e_{ij}$
elements of $M_d(\bc)$, we define {\it Schur multiplication} in
$M_d(\bc)$ as usual,
\begin{equation}\label{schur}
\xb\diamond\yb= \sum_{i,j=1}^d(x_{ij}y_{ij})e_{ij},
\end{equation}
that is, componentwise, $(\xb\diamond \yb)_{ij}:= x_{ij}y_{ij}$.

A linear map $P: M_d(\bc)\to M_d(\bc)$ is said to be {\it Schur
identity-preserving} if its diagonal projection is the identity,
i.e. $\id\diamond P(\id) = \id$ . It is called an {\it entangled}
Markov operator if, in addition, $P(\id)\neq \id$.

The entangled Markov operator (see \cite{AF})  associated to a
stochastic matrix $\Pi=(p_{ij})_{i,j=1}^d$ and to the canonical
systems of matrix units $\{e_{ij}\}_{i,j=1}^d$ of $M_d(\bc)$ is
defined by
\begin{equation}\label{mar1}
P(\xb)_{ij} := \sum_{k,l=1}^d\sqrt{p_{ik}p_{jl}}x_{kl},
\end{equation}
where as before $\xb = \sum_{i,j=1}^d x_{ij}e_{ij}$.

Define a Markov operator $\Psi:M_d(\bc)\to M_d(\bc)$ by
\begin{equation}\label{psi}
\Psi(\xb)=\id\diamond P(\xb), \qquad \xb\in M_d(\bc).
\end{equation}

Given a stochastic matrix $\Pi=(p_{ij})$ put
\begin{equation*}
Fix(\Pi)=\{\p\in\bc^d: \Pi\p=\p\}.
\end{equation*}

To every vector $a=(a_1,\dots,a_d)\in \bc^d$ corresponds a diagonal
matrix $\xb_a$ in $M_d(\bc)$ defined by
\begin{equation}\label{xa}
\xb_a=\left( \begin{array}{cccc} a_1 & 0 & \cdots & 0 \\
0 & a_2 & \cdots & 0 \\
\cdots & \cdots & \cdots & \cdots \\
0 & 0 & \cdots & a_d
\end{array}
\right).
\end{equation}

\begin{lem}\label{fixed} For a Markov operator given by \eqref{psi} one has
$$
M_d(\bc)^\Psi=\{\xb_\p: \ \p\in Fix(\Pi)\}
$$
\end{lem}

\begin{pf} Let $\xb=(x_{ij})\in M_d(\bc)^\Psi$, i.e. $\Psi(\xb)=\xb$. From
\eqref{schur} and \eqref{psi} we conclude that $x_{ij}=0$ if $i\neq
j$. Therefore, due to \eqref{mar1} one finds
\begin{equation*}
\sum_{j=1}^d\sqrt{p_{ij}p_{ij}}x_{jj}=x_{ii}
\end{equation*}
which implies that $(x_{11},\dots,x_{dd})\in Fix(\Pi)$.
\end{pf}

Furthermore, we assume that the dimension of $Fix(\Pi)$ is greater
or equal than 2, i.e. $dim(Fix(\Pi))\geq 2$. Hence, according to
Lemma \ref{fixed} we conclude that $M_d(\bc)^\Psi$ is a non-trivial
commutative subalgebra of $M_d(\bc)$.

\begin{thm}\label{uef}  Let $\Pi$ be a stochastic matrix such that  $dim(Fix(\Pi))\geq 2$.
Then the corresponding Markov operator $\Psi$ given by \eqref{psi}
is uniquely ergodic w.r.t. $M_d(\bc)^\Psi$.
\end{thm}

\begin{pf} To prove the statement, it is enough to establish condition (ii) of
Theorem \ref{unique}. Take any $\xb=(x_{ij})\in M_d(\bc)$. Now we
are going to show that it can be represented as follows
\begin{equation}\label{x12}
\xb=\xb_1+\xb_2,
\end{equation}
where $\xb_1\in M_d(\bc)^\Psi$ and $\xb_2\in \{\yb-\Psi(\yb): \yb\in
M_d(\bc)\}$.

Due to Lemma \ref{fixed} there is a vector $\p\in Fix(\Pi)$ such
that $\xb_1=\xb_\p$, and hence, from \eqref{psi},\eqref{x12} one
finds that
\begin{eqnarray}\label{x-fx1}
\xb_2= \left(
\begin{array}{cccc}
\f_{11} & x_{12} & \cdots & x_{1d} \\
x_{21} & \f_{22} &  \cdots & x_{2d}  \\
\cdots & \cdots & \cdots & \cdots \\
x_{d1} & x_{d2}  & \cdots & \f_{dd}\\
\end{array}
\right),
\end{eqnarray}
where
\begin{equation}\label{x-fx2}
\f_{ii}=\xi_{i}-\sum\limits_{j=1}^d p_{ij}\xi_{j}-\sum\limits_{k,l=1
\atop k\neq j}^d\sqrt{p_{ik}p_{il}}x_{kl}
\end{equation}

The existence of the vectors $\psi=(\psi_1,\dots,\psi_d)$ and
$(\xi_{1},\dots,\xi_{d})$  follows immediately from the following
relations
\begin{equation}\label{x123}
\psi_i+\xi_{i}-\sum\limits_{j=1}^d
p_{ij}\xi_{j}=x_{ii}+\sum\limits_{k,l=1 \atop k\neq
j}^d\sqrt{p_{ik}p_{il}}x_{kl}, \ \ i=1,\dots,d,
\end{equation}
since the number of unknowns is greater than the number of
equations. Note that the equality \eqref{x123} comes from
\eqref{psi}-\eqref{x-fx2}. Hence, one concludes that the equality
\begin{equation*}
M_d(\bc)^{\Psi}+\{\xb-\Psi(\xb):\ \xb\in M_d(\bc)\}=M_d(\bc),
\end{equation*}
which completes the proof.
\end{pf}

Let us provide a more concrete example.

{\bf Example.}  Consider on $M_3(\bc)$ the following stochastic
matrix $\Pi_0$ defined by
\begin{eqnarray}\label{matrix}
\Pi_0=\left(
\begin{array}{ccc}
1 & 0 & 0 \\
0 & 0 & 1\\0 & u & v\\
\end{array}
\right),
\end{eqnarray}
here $u,v\geq 0$, $u+v=1$.

One can immediately find that
\begin{equation}\label{fix}
Fix(\Pi_0)=\{(x,y,y): x,y\in\bc\}.
\end{equation}

Then for the corresponding Markov operator $\Psi_0$, given by
\eqref{psi},\eqref{mar1}, due to Lemma \ref{fixed} one has
\begin{equation}\label{fix2}
M_3(\bc)^{\Psi_0}=\left\{ \left(
\begin{array}{ccc}
x & 0 & 0 \\
0 & y & 0\\
0 & 0 & y\\
\end{array}
\right): x,y\in\bc\right\}.
\end{equation}
So, $M_3(\bc)^{\Psi_0}$ is a nontrivial commutative subalgebra of
$M_3(\bc)$ having dimension 2.

So, according to Theorem \ref{uef} we see that $\Psi_0$ is uniquely
ergodic relative to $M_3(\bc)^{\Psi_0}$. But \eqref{fix2} implies
that $\Psi_0$ is not ergodic. Note that ergodicity of entangled
Markov chains has been studied in \cite{AF}.

\section*{Acknowledgement} The second named author (F.M.) thanks Prof. L. Accardi
for kind hospitality at ``Universit\`{a} di Roma Tor Vergata''
during 18-22 June of 2006. Finally, the authors also would like to
thank to the referee for his useful suggestions which allowed us to
improve the text of the paper.


\begin{thebibliography}{9999}

\bibitem{AH} S.Albeverio, R.H\o egh-Krohn, \textit{Frobenius theory for positive maps of
von Neumann algebras}, Comm. Math. Phys. \textbf{64} (1978),
83--94.

\bibitem{AD} B. Abadie and K. Dykema, \textit{Unique ergodicity of free shifts and some other automorphisms of
$C^*$-algebras},  Jour. Operator Theor. (to appear)
http://www.arxiv.org/math.OA/0608227.

\bibitem{AF} L. Accardi and F. Fidaleo, \textit{Entangled markov chains},  Ann. Mat. Pura Appl.
{\bf 184}(2005), 327--346.


\bibitem{BLRT} D. Berend, M. Lin, J. Rosenblatt and A. Tempelman, \textit{Modulated and
subsequential ergodic theorems in Hilbert and Banach spaces},
Ergod. Th. and Dynam. Sys. \textbf{22} (2002), 1653--1665.

\bibitem{BR} O. Bratteli, D.W. Robinson, \textit{Operator algebras and
quantum statistical mechanics}, I, New York Heidelberg Berlin,
Springer, 1979.

\bibitem{C1} M-D. Choi, \textit{A Schwarz inequality for positive linear maps on
$C^*$-algebras}. {Illinois J. Math.} {\bf 18} (1974), 565--574.

\bibitem{C2}  M-D. Choi, \textit{Completely Positive Linear Maps on Complex Matrices}, {Lin.
Alg. Appl.} {\bf 10}, 285--290 (1975).





\bibitem{FR1} F. Fagnola, R. Rebolledo, \textit{On the existance of stationary states for
quantum dyanamical semigroups}, Jour. Math. Phys. \textbf{42}
(2001), 1296--1308.

\bibitem{FR2} F. Fagnola, R. Rebolledo, \textit{Transience and recurrence of quantum Markov
semi-groups.} Probab. Theory Relat. Fields {\bf 126}(2003),
289-306.

\bibitem{F} F. Fidaleo, \textit{Infinite dimensional entangled Markov
chains},  Random Oper. Stochastic Eq. {\bf 12} (2004), 4, 393--404

\bibitem{FM} F. Fidaleo and F.Mukhamedov, \textit{Strict weak mixing of some C*-dynamical systems based on free
shifts}, Jour. Math. Anal. Appl. {\bf 336}(2007), 180--187.

\bibitem{FV} A. Frigerio, M. Verri, \textit{Long-time asymptotic properties of dynamical semi-groups
on $W^*$-algebras.} Math. Z. {\bf 180}(1982), 275-286.

\bibitem{H} G.H. Hardy, \textit{ Divergent series}, Cambridge Univ. Press, 1949.

\bibitem{I} A. Iwanik, \textit{Unique ergodicity of irreducible Markov operators on $C(X)$},
Studia Math. {\bf 77}(1983), 81--86.

\bibitem{K} V.V. Kozlov, \textit{Weighted means, strict ergodicity and uniform distributions},
Math. Notes, {\bf 78}(2005), 329--337.

\bibitem{KSF} I.P. Kornfeld, Ya.G. Sinai and S.V. Fomin, \textit{ Ergodic Theory},  Springer,
Berlin--Heidelberg--New York, 1982.

\bibitem{LP} R. Longo and C. Peligrad, \textit{Noncommutative topological dynamics and compact actions on
 $C^{*} $-algebras}, J. Funct. Anal. {\bf 58} (1984), 157--174.

\bibitem{MT} F. Mukhamedov and S. Temir, \textit{ A few remarks on mixing properties of
$C^*$-dynamical systems},  Rocky Mount. J. Math. {\bf 37}(2007), 1685--1703.

\bibitem{NSZ} C. Nicolescu, A. Str\"oh, L. Zsid\'o, \textit{Noncommutative
extensions of classical and multiple recurrence theorems,}
J.Operator Theory, {\bf 50}(2003), 3-52.

\bibitem{P} V. Paulsen, \textit{Completely Bounded Maps and Operator Algebras}, Cambridge Univ.
Press, 2002.


\bibitem{T} M. Takesaki, \textit{ Theory of operator algebras I}, Springer,
Berlin-Heidelberg-New York 1979.


\bibitem{W} P. Walters,  \textit{An introduction to ergodic theory},  Springer,
Berlin--Heidelberg--New York, 1982.


\end{thebibliography}
\end{document}